\documentstyle{amsppt} \parindent=24pt
\loadbold \headline={\tenrm\hss\folio\hss}
\magnification=\magstep1\baselineskip=14pt \nopagenumbers\vsize 23 true cm
\hsize 16 true cm \hcorrection{0.9 true cm} \vcorrection{0.4 true cm}

\centerline{\bf ON THE RESIDUALITY A FINITE $\boldkey p$-GROUP} \centerline{\bf
OF $\boldkey H\boldkey N\boldkey N$-EXTENSIONS}
\bigskip
\centerline{D.~I.~Moldavanskii}
\bigskip

\bigskip

\hangindent=28pt \hangafter=-6 {\it A criterion for the $HNN$-extension of a
finite $p$-group to be residually a finite $p$-group is obtained and based on
this criterion the sufficient condition for residuality a finite $p$-group of
$HNN$-extension with arbitrary base group is proved. Then these results are
applied to give for groups from two classes of one-relator groups the necessary
and sufficient condition  to be residually a finite $p$-group.}
\bigskip

\centerline{{\bf 1. Introduction. Statement of results}}
\medskip

Almost all known results on the residual finiteness of generalized free products
of groups are obtained by use of methods offered by G.~Baumslag in work [2].
This methods is based on the assertion (proved in this work)  that the
generalized free product of two finite groups is residual finite and makes use
of the notion of compatible subgroups (introduced there as well). Then this
methods was transferred to the construction of $HNN$-extension of groups: the
residual finiteness of $HNN$-extension of finite group was established
independently in works [1] and [5], and in [1] the sufficient condition of
residual finiteness of $HNN$-extension analogous to corresponding condition from
[2] was given. A criterion of residuality a finite $p$-group of generalized free
product of two finite $p$-groups was obtained by G.~Higman [6] and based on this
criterion certain modification of notion of compatible subgroups leads to
analogous methods of investigation of residuality a finite $p$-group of
generalized free products (see [9]).

The aim of the present paper is to obtain a criterion for $HNN$-extension of a
finite $p$-group to be residually a finite $p$-group (Theorem 1) and to prove
the based on this criterion a sufficient condition (and necessary one) for
residuality a finite $p$-group of any $HNN$-extension (Theorem 2). As
illustration of these results, we give the necessary and sufficient conditions
for groups from two famous classes of one-relator groups to be residually a
finite $p$-group (Theorems 3 and 4).

To formulate results of the paper more explicitly, recall that the chief series
of some group $G$ is a normal series which does not admit a non-trivial normal
supplements. It is easy to see that a normal series of a finite $p$-group is the
chief series if and only if all of its factors are of order $p$. Now, our first
result is

\proclaim{\indent Theorem 1} Let $G$ be a finite $p$-group, let $A$ and $B$ be
subgroups of $G$ and let $\varphi : A\rightarrow B$ be isomorphism. Then
$HNN$-extension $G^{*}=(G,t;\ t^{-1}A\,t=B, \varphi )$ is the residually a
finite $p$-group if and only if there exists a chief series
$$
1=G_{0}\leqslant G_{1}\leqslant \cdots \leqslant G_{n-1}\leqslant G_{n}=G
$$
of group $G$ satisfying the following conditions:
 \roster
\item $(A\cap G_{i})\varphi=B\cap G_{i}$ $(i=0,1,\dots ,n)$;
\item for any $i=0,1,\dots ,n-1$ and for every
$a\in A\cap G_{i+1}$ elements $a\varphi $ and $a$ are congruent modulo subgroup
$G_{i}$.
\endroster
\endproclaim

This Theorem was announced in  [11]. It should be noted that a criterion of
residuality a finite $p$-group of $HNN$-extension of finite $p$-group in another
terms was also obtained in work [12]. Nevertheless, the criterion in Theorem 1
appears to be more suitable for investigation of residuality a finite $p$-group
of $HNN$-extensions with infinite base group (see Theorem 2 below). Let us
remark also that the proof of Theorem 1 is quite elementary since it makes use
of only ordinary properties of construction of $HNN$-extension. By the similar
arguments one can also prove Higman's Theorem mentioned above. It is relevant to
recall that original proof of Higman, as well as the proof of respective result
in [12], exploits the construction of wreath product.
\smallskip

Higman's Theorem implies, in particular, that the generalized free product of
two finite $p$-groups is residually a finite $p$-group provided that amalgamated
subgroups are cyclic. The following simple example demonstrates that, unlike of
this, the assumption that subgroups $A$ and $B$ of a finite $p$-group $G$ are
cyclic does not guarantee the residuality a finite $p$-group of group
$G^{*}=(G,t;\ t^{-1}A\,t=B, \varphi )$.

Let $H^{*}=(H,t;\ t^{-1}a\,t=b^{p})$ be the $HNN$-extension of group
$$
H=\langle a,b;\ a^{-1}ba=b^{1+p},\ a^{p}=1,\ b^{p^{2}}=1\rangle
$$
of order $p^{3}$, where associated subgroups $A$ and $B$ are generated by
elements $a$ and $b^{p}$ respectively. If we suppose that the group $H^{*}$ is
residually a finite $p$-group then intersection of all members
$\gamma_{n}(H^{*})$ of its lower central series must coincide with identity
subgroup. Therefore, we can choose the number $n$ such that $a\in
\gamma_{n}(H^{*})\setminus \gamma_{n+1}(H^{*})$. But then congruences $b\equiv
b^{1+p}\pmod {\gamma_{n+1}(H^{*})}$ and $a\equiv b^{p}\pmod
{\gamma_{n+1}(H^{*})}$ will imply that $a\in \gamma_{n+1}(H^{*})$. (Conditions
of Theorem 1 are not satisfied here since $A\cap B=1$ and the first non-identity
member of any chief series of group $H$ must coincide with its centre $B$.)

Nevertheless, it is relevant to expect that in the case when subgroups $A$ and
$B$ are cyclic one can find a more simple criterion for group $G^{*}$ to be
residually a finite $p$-group. As some confirmation of that, in the case when
subgroups $A$ and $B$ are equal we have

\proclaim{\indent Corollary} Let $G$ be a finite $p$-group, let $A$ be
non-identity cyclic subgroup of $G$ with generator $a$ and let $\varphi $ be
automorphism of $A$ such that $a\varphi =a^{k}$ for some integer $k$ (coprime to
$p$). Then group $G^{*}=(G,t;\ t^{-1}a\,t=a^{k})$ is residually a finite
$p$-group if and only if $k\equiv 1\pmod p$.
\endproclaim

Indeed, if $1=G_{0}\leqslant G_{1}\leqslant \cdots \leqslant G_{n-1}\leqslant
G_{n}=G$ is a chief series of group $G$ then distinct members of sequence
$(A\cap G_{i})$ $(i=0,1,\dots ,n)$ constitute the unique chief series of group
$A$. Hence $(A\cap G_{i})\varphi=A\cap G_{i}$ for all $i=0,1,\dots ,n$. If $a\in
G_{i+1}\setminus G_{i}$ then, since element $aG_i$ of quotient group $G/G_i$ is
of order $p$, the equality $(a\varphi )G_{i}=aG_{i}$ implies that $k\equiv
1\pmod p$. Conversely, if this congruence is fulfilled then it is obvious that
any chief series of group $G$ satisfies the condition (2) of Theorem 1.
\smallskip

In what follows we need some notions that came from the work [2] and are
utilized now in almost all investigations of residual properties of
$HNN$-extensions.

A family $(N_{\lambda})_{\lambda \in \Lambda}$ of normal subgroups of a group
$G$ is said to be filtration if $\bigcap _{\lambda \in \Lambda}N_{\lambda}=1$.
If $H$ is a subgroup of group $G$ and if $\bigcap _{\lambda \in
\Lambda}HN_{\lambda}=H$ then filtration $(N_{\lambda})_{\lambda \in \Lambda}$ is
called $H$-filtration. If $H$ and $K$ are two subgroups of $G$ then filtration
will be called $(H,K)$-filtration if it is $H$-filtration and $K$-filtration
simultaneously.

Let now $G$ be a group with subgroups $A$ and $B$ and let  $\varphi :
A\rightarrow B$ be isomorphism. Subgroup $H$ of group $G$ is called
$(A,B,\varphi )$-compatible if $(A\cap H)\varphi =B\cap H$. (So, the condition
(1) in Theorem 1 means that all subgroups $G_i$ are $(A,B,\varphi
)$-compatible.) It is easy to see that if $H$ is a normal $(A,B,\varphi
)$-compatible subgroup of $G$ then the mapping $\varphi_{{}_H} : AH/H\to BH/H$
(well) defined by the rule $(aH)\varphi_{{}_{H}}=(a\varphi )H$ (where $a\in A$)
is an isomorphism of subgroup $AH/H$ of quotient group $G/H$ onto subgroup
$BH/H$. Furthermore, the natural homomorphism of group $G$ onto quotient group
$G/H$ can be extended to homomorphism $\rho_{{}_H}$ (sending $t$ to $t$) of
$HNN$-extension $G^{*}=(G,t;\ t^{-1}A\,t=B, \varphi )$ onto $HNN$-extension
$$
G^{*}_{H}=(G/H,t;\ t^{-1}AH/H\,t=BH/H,\, \varphi_{{}_H}).
$$

Let $\Cal F_{G}(A,B,\varphi )$ denote the family of all $(A,B,\varphi
)$-compatible normal subgroups of finite index of group $G$. Mentioned above
sufficient condition of residual finiteness of $HNN$-extension $G^{*}$ of group
$G$ in [1] consists of requirement for the family $\Cal F_{G}(A,B,\varphi )$ to
be $(A,B)$-filtration. To formulate the analogous condition for residuality a
finite $p$-group of $HNN$-extension we give here  based on the Theorem 1
corresponding modification of notion $(A,B,\varphi )$-compatibility.

Let, as before, $G$ be a group with subgroups $A$ and $B$ and $\varphi :
A\rightarrow B$ be isomorphism. Let $p$ be a prime number. Subgroup $H$ of group
$G$ will be called $(A,B,\varphi,p)$-compatible if there exists a sequence
$$
H=G_{0}\leqslant G_{1}\leqslant \cdots \leqslant G_{n-1}\leqslant G_{n}=G
$$
of subgroups of group $G$ such that

\hangindent=28pt \hangafter=-6\noindent 1) for any $i=0,1,\dots ,n$ subgroup
$G_{i}$ is $(A,B,\varphi )$-compatible and normal in $G$, and

\hangindent=28pt \hangafter=-6\noindent 2) for any $i=0,1,\dots ,n-1$ the
quotient group $G_{i+1}/G_{i}$ is of order $p$ and for every $a\in A\cap
G_{i+1}$ elements $a\varphi $ and $a$ are congruent modulo subgroup $G_{i}$.
\smallskip

Let $\Cal F_{G}^{p}(A,B,\varphi )$ denote the family of all $(A,B,\varphi
,p)$-compatible subgroups of group $G$. Thus, Theorem 1 asserts, actually, that
if $G$ is a finite $p$-group then $HNN$-extension $G^{*}=(G,t;\ t^{-1}A\,t=B,
\varphi )$ is residually a finite $p$-group if and only if the family $\Cal
F_{G}^{p}(A,B,\varphi )$ contains identity subgroup.

The following Theorem giving sufficient condition for $HNN$-extension to be
residually a finite $p$-group can be considered, as well, as a confirmation of
the fact that the notion of $(A,B,\varphi ,p)$-compatibility indeed can be used
as $p$-analog of notion of $(A,B,\varphi )$-compatibility.

\proclaim{\indent Theorem 2} Let $G$ be a group with subgroups $A$ and $B$ and
let $\varphi : A\rightarrow B$ be isomorphism. Let $G^{*}=(G,t;\ t^{-1}A\,t=B,
\varphi )$ be the $HNN$-extension of $G$. Then \roster
\item if group $G^{*}$ is residually a finite $p$-group then the family
$\Cal F_{G}^{p}(A,B,\varphi)$ is a filtration;

\item if the family $\Cal F_{G}^{p}(A,B,\varphi )$ is an $(A,B)$-filtration then
group $G^{*}$ is residually a finite $p$-group.
\endroster
\endproclaim

In the case when group $G$ is abelian and $A$ and $B$ are proper subgroups of
$G$ one can say somewhat more. Let $g\in G\setminus A$ and $h\in G\setminus B$.
Then commutator $u=[t^{-1}gt, h]$ is not equal to identity since its expression
$u=t^{-1}g^{-1}t\,h^{-1}\,t^{-1}gt\,h$  is reduced in group $G^{*}$. If the
group $G^{*}$ is residually a finite $p$-group then some normal subgroup $N$ of
finite $p$-index of $G^{*}$ does not contain element $u$. If $M=G\cap N$ then
since quotient group $G/M$ is abelian it follows that $g\notin AM$. Thus, since
subgroup $M$ is $(A,B,\varphi ,p)$-compatible (see Lemma 2.2 below), we have
obtained the

\proclaim{\indent Corollary} If group $G$ is abelian and $A$ and $B$ are proper
subgroups of $G$ then group $G^{*}$ is residually a finite $p$-group if and only
if the family $\Cal F_{G}^{p}(A,B,\varphi )$ is $(A,B)$-filtration.
\endproclaim

Let us consider now two classes of one-relator groups. The first is the class of
Baumslag -- Solitar groups, i.~e. class of groups with presentation
$$G(l,m)=\langle a, b;\ b^{-1}a^{l}b=a^{m}\rangle,$$ where without loss of
generality one can assume that $|m|\geqslant l>0$. It is well-known (see [3,
10]) that the group $G(l,m)$ is residually finite if and only if either $l=1$,
or $|m|=l$.

The second class consists of certain $HNN$-extensions of Baumslag -- Solitar
groups, namely, of groups with presentation
$$
\multline
H(l,m;\,k)=\langle a,t;\ t^{-1}a^{-k}t\,a^{l}\,t^{-1}a^{k}t=a^{m}\rangle \\
=\langle a,b,t;\ b^{-1}a^{l}\,b=a^{m},\ t^{-1}a^{k}\,t=b\rangle,
\endmultline
$$
where (again without loss of generality) it can be supposed that $|m|\geqslant
l>0$ and $k>0$. Some properties of these groups were established by
A.~M.~Brunner [4] (see also [8]). In particular, it is known that the group
$H(l,m;\,k)$ is residually finite if and only if $|m|=l$.

Theorems 1 and 2 will be applied here to prove following assertions:

\proclaim{\indent Theorem 3} Group $G(l,m)=\langle a, b;\
b^{-1}a^{l}b=a^{m}\rangle$ (where $|m|\geqslant l>0$) is residually a finite
$p$-group if and only if either $l=1$ and $m\equiv 1\pmod p$, or $|m|=l=p^{r}$
for some $r\geqslant 0$ and also if $m=-l$ then $p=2$.
\endproclaim

\proclaim{\indent Theorem 4} Group $H(l,m;\,k)=\langle a,t;\
t^{-1}a^{-k}t\,a^{l}\,t^{-1}a^{k}t=a^{m}\rangle$ (where $|m|\geqslant l>0$ and
$k>0$) is residually a finite $p$-group if and only if $|m|=l=p^{r}$ and
$k=p^{s}$ for some integers $r\geqslant 0$ and $s\geqslant 0$ and also if $m=-l$
then $p=2$ and $s\leqslant r$.
\endproclaim
\bigskip

\centerline{{\bf 2. Proof of Theorems 1 and 2}}
\medskip

To prove Theorem 1 we begin with simple and well-known (see, e.~g., [12,
Proposition 1] remark:

\proclaim{\indent Lemma 2.1} Let $G$ be a finite $p$-group, $A,B\leqslant G$ and
let $\varphi : A\rightarrow B$ be isomorphism. Group $G^{*}=(G,t;\ t^{-1}A\,t=B,
\varphi )$ is residually a finite $p$-group if and only if there exists a
homomorphism $\rho$ of group $G^{*}$ onto some finite $p$-group $X$ such that
$\text{\rm Ker}\,\rho\cap G=1$.
\endproclaim

In fact, the part \lq\lq only if\rq\rq\ of Lemma is obvious since group $G$ is
finite. Conversely, if $\text{Ker}\,\rho\cap G=1$ then (see [7]) subgroup
$\text{Ker}\,\rho$ is free. So, group $G^{*}$ is free-by-(finite $p$-group) and
therefore is residually a finite $p$-group.

Suppose now that the $HNN$-extension $G^{*}=(G,t;\ t^{-1}A\,t=B,\varphi )$ of
finite $p$-group $G$ is residually a finite $p$-group. Then by Lemma 2.1 we can
consider the group $G$ as a subgroup of some finite $p$-group $X$ with element
$x$ such that $x^{-1}ax=a\varphi$ for all $a\in A$. Let
$$
1=X_{0}\leqslant X_{1}\leqslant \cdots \leqslant X_{n-1}\leqslant X_{n}=X
$$
be a chief series of group $X$ and $G_{i}=G\cap X_{i}$ ($i=0,1,\dots , n$). Then
distinct members of sequence $G_{0}$, $G_{1}$, \dots , $G_{n-1}$, $G_{n}$ form
the chief series of group $G$. Since $A\cap G_{i}=A\cap X_{i}$ and $B\cap
G_{i}=B\cap X_{i}$,
$$
(A\cap G_{i})\varphi=(A\cap X_{i})\varphi =(A\cap X_{i})^{x}=A^{x}\cap
X_{i}=B\cap X_{i}=B\cap G_{i}.
$$

Let $\psi $ be the embedding of quotient group $G/G_{i}$ into quotient group
$X/X_{i}$ which takes coset $gG_{i}$ to coset $gX_{i}$. Since subgroup
$(G_{i+1}/G_{i})\psi $ is contained in central subgroup $X_{i+1}/X_{i}$ of
$X/X_{i}$, for any element $a\in A\cap G_{i+1}$ we have
$$
(aG_{i})\psi =aX_{i}=a^{x}X_{i}=(a\varphi)X_{i}=((a\varphi)G_{i})\psi.
$$
Since the mapping $\psi $ is injective, this implies that
$(a\varphi)G_{i}=aG_{i}$. Thus, we have the chief series of group $G$ satisfying
the conditions (1) and (2) of Theorem 1.
\smallskip

Conversely, suppose that some chief series
$$
1=G_{0}\leqslant G_{1}\leqslant \cdots \leqslant G_{n-1}\leqslant G_{n}=G
$$
of group $G$ satisfies the conditions (1) and (2) of Theorem 1. By induction on
$n$ we shall show that then there exists a homomorphism of group $G^{*}$ onto
some finite $p$-group $X$ which acts injectively on subgroup $G$ (and thus, in
view of Lemma 2.1, we shall prove that the group $G^{*}$ is residually a finite
$p$-group).

It is easy to see that if $n=1$ then for group $X$ we can take the direct
product of group $G$ and cyclic group of order $p$ with generator $x$; the
desired mapping acts on group $G$ identically and sends element $t$ to element
$x$.

Let $n>1$. Since $(A\cap G_{1})\varphi=B\cap G_{1}$ and subgroup $G_{1}$ is of
order $p$, we must consider only two following cases:

\noindent a) $G_{1}\leqslant A$ ¨ $G_{1}\leqslant B$;

\noindent b) $A\cap G_{1}=B\cap G_{1}=1$.

In the case a) we set $\overline G=G/G_{1}$, $\overline G_{i}=G_{i}/G_{1}$
$(i=1,2,\dots , n)$, $\overline A=A/G_{1}$ and $\overline B=B/G_{1}$. Then
$$1=\overline G_{1}\leqslant \overline G_{2}\leqslant \cdots \leqslant \overline
G_{n-1}\leqslant \overline G_{n}=\overline G$$ is the chief series of group
$\overline G$. Since subgroup $G_{1}$ is $(A, B, \varphi )$-compatible, the
mapping $\overline \varphi =\varphi _{{}_{G_{1}}}$ is an isomorphism of subgroup
$\overline A$ onto subgroup $\overline B$. Moreover, since $\overline A\cap
\overline G_{i}=(A\cap G_{i})/G_{1}$ ¨ $\overline B\cap \overline G_{i}=(B\cap
G_{i})/G_{1}$, we have $(\overline A\cap \overline G_{i})\overline \varphi=
\overline B\cap \overline G_{i}$. Also, it can be immediately checked that for
any $i=1,2,\dots , n-1$ and for every element $aG_{1}\in \overline A\cap
\overline G_{i+1}$ cosets $aG_{1}$ and $(aG_{1})\overline \varphi$ are congruent
modulo subgroup $\overline G_{i}$. Consequently, by induction there exists a
homomorphism  $\sigma $ of group $\overline G^{*}= (\overline G,t;\
t^{-1}\overline A\,t=\overline B, \overline \varphi )$ onto finite $p$-group $X$
such that $\text{Ker}\,\sigma\cap \overline G=1$.

Let $\rho  =\rho _{{}_{G_{1}}}$ be the homomorphism of group $G^{*}$ onto the
group $\overline G^{*}$ which extends natural mapping of $G$ onto quotient group
$\overline G$. Let $L=\text{Ker}\,(\rho \sigma)$. Then $\text{Ker}\,\rho
=G_{1}$, $G^{*}/L\simeq X$ ¨ $G\cap L=G_{1}$. Furthermore, since $L/G_{1}\simeq
\text{Ker}\,\sigma$ and the group $\text{Ker}\,\sigma$ is free there exists a
free subgroup $F$ of $L$ such that $L=FG_{1}$ and $F\cap G_{1}=1$. As the
condition (2) implies that subgroup $G_{1}$ belongs to the centre of $G^{*}$, we
have $L=F\times G_{1}$. Let $N$ denote intersection of all normal subgroups of
index $p$ of group $L$. Then $N$ is contained in $F$ and is normal in $G$
because it is characteristic in $L$. Moreover, quotient group $L/N$ is a finite
$p$-group since group $L$ is finitely generated as a subgroup of finite index of
finitely generated group $G^{*}$. At last, $N\cap G=N\cap L\cap G=N\cap
G_{1}=N\cap F\cap G_{1}=1$. Thus, the natural homomorphism of group $G^{*}$ onto
quotient group $G^{*}/N$ is desired.

In the case b) we set $A_{1}=AG_{1}$ and $B_{1}=BG_{1}$. Since $G_{1}$ is a
central subgroup of $G$, $A_{1}=A\times G_{1}$ ¨ $B_{1}=B\times G_{1}$.
Therefore, the mapping $\varphi_{1} :A_{1}\rightarrow B_{1}$ which carries
element $g\in A_{1}$, $g=ax$ (where $a\in A$ and $x\in G_{1}$), to element
$(a\varphi )x$ is an isomorphism. Since $A_{1}\cap G_{i}=(A\cap G_{i})G_{1}$ and
$B_{1}\cap G_{i}=(B\cap G_{i})G_{1}$, we have $(A_{1}\cap
G_{i})\varphi_{1}=B_{1}\cap G_{i}$. Furthermore, if element $g=ax$ (where $a\in
A$ ¨ $x\in G_{1}$) belongs to subgroup $A_{1}\cap G_{i+1}=(A\cap G_{i+1})G_{1}$
then $a\in A\cap G_{i+1}$ and, therefore, $(g\varphi_{1})G_{i}=(a\varphi
)xG_{i}=(a\varphi )G_{i}\cdot xG_{i} =aG_{i}\cdot xG_{i}=gG_{i}$. Hence, by the
case a), treated above,  there exists a homomorphism $\sigma$ of group
$$
G^{*}_{1}=(G,t;\ t^{-1}A_{1}t=B_{1}, \varphi_{1})
$$
onto some finite $p$-group $X$ which acts injectively on subgroup $G$. Since the
isomorphism $\varphi$ coincides with restriction to subgroup $A$ of isomorphism
$\varphi_{1}$, the identity mapping of group $G$ can be extended to homomorphism
$\rho: G^{*}\rightarrow G^{*}_{1}$. Then the homomorphism $\rho \sigma$ maps
group $G^{*}$ onto group $X$ and acts injectively on $G$. Thus, inductive step
is completed and Theorem 1 is proved.
\smallskip

Let us proceed to prove Theorem 2. Let $G$ be a group with subgroups  $A$ and
$B$ and let $\varphi:A\rightarrow B$ be isomorphism.

\proclaim{\indent Lemma 2.2} {\rm a)} Any normal $(A,B,\varphi )$-compatible
subgroup $H$ of group $G$ belongs to the family $\Cal F_{G}^{p}(A,B,\varphi )$
if and only if the group $G^{*}_{H}$ is residually a finite $p$-group.

{\rm b)} Let $N$ be a normal subgroup of finite $p$-index of group $G^{*}=(G,t;\
t^{-1}A\,t=B, \varphi )$ and $M=G\cap N$. Then $M\in \Cal F_{G}^{p}(A,B,\varphi
)$.

{\rm c)} The family $\Cal F_{G}^{p}(A,B,\varphi )$ is closed under intersections
of finite collections of subgroups.
\endproclaim

The proof of all assertions in Lemma 2.2 does not evoke a special difficulties.
The validity of item  ) follows immediately from Theorem 1 applied to group
$G/H$ and from definition of $(A,B,\varphi,p)$-compatibility. As well, item b)
is direct consequence of item a) and of Lemma 2.1. To prove item c) it is
sufficient to remark that if subgroups $H$ and $K$ belong to family $\Cal
F_{G}^{p}(A,B,\varphi )$ and $L=H\cap K$ then there exists a homomorphism $\rho$
of group $G^*_L$ into direct product $G^*_H\times G^*_K$ which acts on subgroup
$G/L$ injectively.

The assertion (1) in Theorem 2 is an obvious consequence of item b) in Lemma
2.2. It follows from the item  ) in Lemma 2.2  that to prove the assertion (2)
it is enough to show that for any non-identity element $w$ of group $G^*$ there
exists subgroup $H\in \Cal F_{G}^{p}(A,B,\varphi )$ such that $w\rho_{{}_H}\neq
1$.

If element $w$ belongs to subgroup $G$ then the existence of such subgroup is
ensured by assumption that the family $\Cal F_{G}^{p}(A,B,\varphi )$ is
filtration.

 Let $w=g_0t^{\varepsilon_1}g_1t^{\varepsilon_2}g_2\cdots
t^{\varepsilon_n}g_n$ be a reduced form of element $w$, where $n\geqslant 1$.
Then for any $i=1, 2, \dots , n-1$ such that $\varepsilon_i + \varepsilon_{i+1}
= 0$ we have $g_i\notin A$ if $\varepsilon_i = -1$ and $g_i\notin B$ if
$\varepsilon_i =1$. Now, from the assumption that the family $\Cal
F_{G}^{p}(A,B,\varphi )$ is $(A,B)$-filtration and from item c) in Lemma 2.2 it
follows that there exists subgroup $H\in \Cal F_{G}^{p}(A,B,\varphi )$ such that
for any $i=1, 2, \dots , n-1$ $g_i\notin AH$ if $g_i\notin A$ and $g_i\notin BH$
if $g_i\notin B$. This means that the expression
$(g_0H)t^{\varepsilon_1}(g_1H)t^{\varepsilon_2}(g_2H)\cdots
t^{\varepsilon_n}(g_nH)$ is a reduced form of element $w\rho_{{}_H}$ and
therefore this element is not equal to 1. The proof of Theorem 2 is complete.
\bigskip

\centerline{{\bf 3. Proof of Theorems 3 ¨ 4}}
\medskip

The group $G(l,m)=\langle a, b;\ b^{-1}a^{l}b=a^{m}\rangle$ is an
$HNN$-extension of infinite cyclic group $A$ generated by element $a$ with
stable letter $b$, associated subgroups $A^l$ and $A^m$ and with isomorphism
$\varphi$ which sends element $a^l$ to element $a^m$. If this group is
residually a finite $p$-group then it is residually finite and therefore, as was
remarked above,  $l=1$ or $|m|=l$ (recall that we assume that $|m|\geqslant
l>0$).

Suppose, at first, that $l=1$. Let $\rho $ be the homomorphism of group $G(1,m)$
onto finite $p$-group $X$ such that $a\rho \neq 1$. If $p^{s}$ denotes the order
of element $a\rho $ then $\rho $ passes through the group
$$
G_{s}(1,m)=\langle a, b;\ b^{-1}ab=a^{m},\ a^{p^{s}}=1\rangle
$$
which is an $HNN$-extension of finite cyclic group $A/A^{p^s}$. Since there
exists a homomorphism of group $G_{s}(1,m)$ onto finite $p$-group $X$ which acts
injectively on base group of this $HNN$-extension, by Lemma 2.1 group
$G_{s}(1,m)$ is residually a finite $p$-group. The Corollary from Theorem 1
implies now that $m\equiv 1\pmod p$.

Conversely, if the congruence $m\equiv 1\pmod p$ holds and therefore (by the
same Corollary) any group $G_{s}(1,m)$ is residually a finite $p$-group then the
group $G(1,m)$ is residually a finite $p$-group too, since it, as easy to see,
is residually groups $G_{s}(1,m)$ $(s=1,2,\dots\,)$.

If $l>1$ (and $|m|=l$) then by Corollary from Theorem 2 group $G(l,m)$ is
residually a finite $p$-group if and only if the family $\Cal
F_{A}^{p}(A^{l},A^{m},\varphi )$ is $(A^{l},A^{m})$-filtration.

Let $l=l_1p^r$ where $r\geqslant 0$ and $(l_1, p)=1)$. It is evident that if
subgroup $A^k$ of group $A$ is $(A^{l},A^{m},p)$-compatible then $k$ is a
$p$-number. Also, it is easy to see that if $m=l$ then any subgroup of form
$A^{p^s}$ is $(A^{l},A^{m},p)$-compatible and if $m=-l$ and $s>r$ then subgroup
$A^{p^s}$ is $(A^{l},A^{m},p)$-compatible if and only if $p=2$. Let us note, at
last, that if integers $x$ and $y$ are such that $l_1x+p^sy=1$ then the equality
$a^{p^r}=(a^l)^x\cdot (a^{p^s})^{p^ry}$ holds and therefore $a^{p^r}\in A^l\cdot
A^{p^s}$. Thus, the family $\Cal F_{A}^{p}(A^{l},A^{m},\varphi )$ is
$(A^{l},A^{m})$-filtration if and only if $l=p^{r}$ for some $r\geqslant 0$ and
also if $m=-l$ then $p=2$. Theorem 3 is proved.
\smallskip

Let us proceed now to prove the Theorem 4. Suppose that the group
$$H(l,m;\,k)=\langle a,b,t;\ b^{-1}a^{l}b=a^{m},\ t^{-1}a^{k}t=b\rangle$$ (where
$|m|\geqslant l>0$ ¨ $k>0$) is residually a finite $p$-group. Then (see [1, 7])
$|m|=l$ and since group $G(l,m)$ is a subgroup of $H(l,m;\,k)$ the Theorem 3
implies that $|m|=l=p^{r}$ for some $r\geqslant 0$ and also if $m=-l$ then
$p=2$.

Let $k=p^{s}k_{1}$, where $s\geqslant 0$ and $(k_{1},p)=1$. If $k_{1}>1$ then a
$t$-reduced form of element $w=\bigl[t^{-1}a^{-p^{s}}t\,a^{l}t^{-1}a^{p^{s}}t,
a\bigr]$ of group $H(l,m;\,k)$ is of length 8 and therefore $w\neq 1$. On the
other hand, if we suppose that $N$ is a normal subgroup of group $H(l,m;\,k)$
such that $a^{p^n}\in N$ for some $n\geqslant 0$, then $a^{p^{s}}\equiv
a^{kc}\pmod N$, where $c$ is an integer satisfying the congruence $k_{1}c\equiv
1\pmod {p^{n}}$. Hence $t^{-1}a^{p^{s}}t\equiv b^{c}\pmod N$, and since
$m=\varepsilon l$ for some $\varepsilon=\pm 1$, we have
$t^{-1}a^{-p^{s}}t\,a^{l}t^{-1}a^{p^{s}}t\equiv a^{\varepsilon ^{c}l} \pmod N$
and therefore $w\in N$. Consequently, $w$ belongs to every normal subgroup of
finite $p$-index of group $H(l,m;\,k)$, but this is impossible since this group
is residually a finite $p$-group. Thus, $k=p^{s}$.

Suppose that $m=-l$ (and therefore $p=2$). Let $\sigma $ be homomorphism of
group $H(2^{r},-2^{r};\,2^{s})$ onto finite $p$-group $X$ such that $b\sigma
\neq 1$. Let
$$1=X_{0}\leqslant X_{1}\leqslant \cdots \leqslant X_{n}=X$$
be a
chief series of group $X$ and $b\sigma \in X_{i+1}\setminus X_{i}$. Since
$X_{i+1}/X_{i}$ is a central subgroup of group $X/X_{i}$ we have
$(a\sigma)^{2^{r}}\equiv (a\sigma)^{-2^{r}}\pmod {X_{i}}$ and
$(a\sigma)^{2^{s}}\equiv b\sigma \pmod {X_{i}}$. Therefore, element
$(a\sigma)^{2^{r+1}}$ belongs to subgroup $X_{i}$ and element
$(a\sigma)^{2^{s}}$ does not belong to that subgroup. This implies the
inequality $s\leqslant r$.

Conversely, let us show that for any prime number $p$ and for any integers
$r\geqslant 0$ and $s\geqslant 0$ the group $H(p^{r},\varepsilon p^{r};\,p^{s})$
(where $\varepsilon =\pm 1$ and if $\varepsilon =-1$ then $p=2$ and $s\leqslant
r$) is residually a finite $p$-group.

Let $G=\langle a, b;\ b^{-1}a^{p^{r}}b=a^{\varepsilon p^{r}}\rangle$ and let $A$
and $B$ be cyclic subgroups of $G$ generated by elements $a$ and $b$
respectively and $A_{1}=A^{p^{s}}$. Then the group $H(p^{r},\varepsilon
p^{r};\,p^{s})$ is $HNN$-extension $(G,t;\,t^{-1}A_{1}t=B,\,\varphi)$ where
isomorphism $\varphi $ is defined by the equality $a^{p^{s}}\varphi =b$.

By Theorem 2 it is enough to show that the family $\Cal
F_{G}^{p}(A_{1},B,\varphi )$ is $(A_{1},B)$-filtration. Nevertheless, we begin
with somewhat weaker assertion.

\proclaim{\indent Lemma 3.1} The family of all normal $(A_{1},B)$-compatible
subgroups of finite $p$-index of group $G$ is $(A_{1},B)$-filtration.
\endproclaim

\demo{\indent Proof} It is evident that the group $G$ is amalgamated free
product $$G=(A*K;\ a^{p^{r}}=c)$$ of groups $A$ and $K=\langle b,c;\
b^{-1}cb=c^{\varepsilon}\rangle$. For any integer $n>\text{max}\,(r,s)$ let
$L_{n}$ be subgroup of $K$ generated by elements $b^{p^{n-s}}$ and
$c^{p^{n-r}}$. Since when $\varepsilon =-1$ then $p=2$, in any case $L_{n}$ is a
normal subgroup of $K$. Also,  it is obvious that the quotient group $K/L_{n}$
is of order $p^{2n-r-s}$ and its elements $bL_{n}$ and $cL_{n}$ are of order
$p^{n-s}$ and $p^{n-r}$ respectively. Therefore, we can construct amalgamated
free product $G_{n}=(A/A^{p^{n}}*K/L_{n};\ (aA^{p^{n}})^{p^{r}}=cL_{n})$.

Let $\rho_{n}$ be homomorphism of group $G$ onto group $G_{n}$ which extends the
natural mappings of group $A$ onto quotient $A/A^{p^{n}}$ and of group $K$ onto
quotient $K/L_{n}$. We claim that for any element $g\in G$ if $g\neq 1$ or if
$g\notin A_{1}$ or if $g\notin B$ then there exists an integer $n$ such that
$g\rho_{n}\neq 1$ or $g\rho_{n}\notin A_{1}\rho_{n}$ or $g\rho_{n}\notin
B\rho_{n}$ respectively. In fact, this is evident if $g\in A$ or $g\in K$.
Moreover, in these cases the following assertion is obvious too: if element $g$
does not belong to amalgamated subgroup in above decomposition of group $G$,
then for all numbers $n$ that are large enough element $g\rho_{n}$ does not
belong to amalgamated subgroup in decomposition of group $G_{n}$. This implies
that if reduced form of element $g$ is of length $>1$ then for suitable $n$ the
reduced form of element $g\rho_{n}$ is of the same length. Therefore, this
element does not belong to both subgroups $A_{1}\rho_{n}$ and $B\rho_{n}$.

Let $g$ be any non-identity element of group $G$ and let integer $n$ be chosen
so that $g\rho_{n}\neq 1$. Since the group $G_{n}$ (by Higman's criterion
mentioned above) is residually a finite $p$-group and its subgroups
$A/A^{p^{n}}$ and $K/L_{n}$ are finite, there exists a normal subgroup $M$ of
finite $p$-index of group $G_{n}$ such that $g\rho_{n}\notin M$ and
$A/A^{p^{n}}\cap M = K/L_{n}\cap M=1$. Let $N=M\rho_{n}^{-1}$. Then $N$ does not
contain element $g$ and is a normal $(A_{1},B,\varphi )$-compatible subgroup of
finite $p$-index of group $G$. The similar arguments show that if element $g$
does not belong to subgroup $A_{1}$ or to subgroup $B$ then subgroup $N$ can be
chosen so that element $g$ does not belong to subgroup $A_{1}N$ or to subgroup
$BN$ respectively.
\enddemo

\proclaim{\indent Lemma 3.2} For any integer $n>s$ there exists a normal
subgroup $N$ of finite $p$-index of group $G$ such that $A\cap N=A^{p^{n}}$,
$B\cap N=B^{p^{n-s}}$ and $a^{p^{n-1}}\equiv b^{p^{n-s-1}}\pmod N$.
\endproclaim

\demo{\indent Proof} Firstly suppose that $n>r$. Assuming all notations from the
proof of Lemma 3.1 let us consider also cyclic subgroup $D$ of group $K$,
generated by element $d=c^{p^{n-r-1}}b^{-p^{n-s-1}}$. If in the case
$\varepsilon =-1$ we shall suppose that $n-s-1>0$ then in any case
$DL_{n}/L_{n}$ will be a central subgroup of order $p$ of group $K/L_{n}$ and
elements $c(DL_{n})$ and $b(DL_{n})$ of quotient group $K/DL_{n}$ will be of
order $p^{n-r}$ and $p^{n-s}$ respectively. Since the group $G_{n}^{\prime
}=(A/A^{p^{n}}*K/DL_{n};\, (aA^{p^{n}})^{p^{r}}=c(DL_{n}))$ is residually a
finite $p$-group, it contains a normal subgroup $M$ of finite $p$-index such
that $A/A^{p^{n}}\cap M=K/DL_{n}\cap M=1$. Then the preimage $N$ of subgroup $M$
under the obvious homomorphism of group $G$ onto group $G_{n}^{\prime }$ is
desired subgroup.

If $\varepsilon =-1$ and $n-s-1=0$ then since $n>r$ and $s\leqslant r$ we have
$s=r$ and $n=r+1$. In this case desired subgroup $N$ coincides with kernel of
homomorphism $\sigma $ of group $G$ onto cyclic group $X$ of order $2^{r+1}$
with generator $x$, where $a\sigma =x$ and $b\sigma =x^{2^{r}}$.

At last, if $n\leqslant r$ then group $G$ can be mapped onto group
$$
T=\langle a,b;\ a^{p^{n}}=1,\ b^{p^{n-s}}=1, \ a^{p^{n-1}}=b^{p^{n-s-1}}\rangle.
$$
Since group $T$ is residually a finite $p$-group it has a normal subgroup of
finite $p$-index whose intersections with (finite) free factors are trivial. The
preimage of this subgroup is desired subgroup of $G$.
\enddemo

Now we are able to prove the assertion mentioned above and thus to complete the
proof of Theorem 4.

\proclaim{\indent Lemma 3.3} The family $\Cal F_{G}^{p}(A_{1},B,\varphi )$ is
$(A_{1},B)$-filtration.
\endproclaim

\demo{\indent Proof} We shall show that any $(A_{1},B)$-compatible subgroup $N$
of finite $p$-index of group $G$ such that $n>s$ where the integer $n$ is
determined by the equality $A\cap N=A^{p^{n}}$ contains some subgroup $M$ from
family $\Cal F_{G}^{p}(A_{1},B,\varphi )$. It is evident that then the desired
assertion will be applied by Lemma 3.1.

Since inequalities $k<n-s$ and $n-k>s$ are equivalent, it follows from Lemma 3.2
that for every integer $k$, $0\leqslant k<n-s$, in group $G$ there exists a
normal subgroup $N_{k}$ of finite $p$-index such that $A\cap N_{k}=A^{p^{n-k}}$,
$B\cap N_{k}=B^{p^{n-k-s}}$ ¨ $a^{p^{n-k-1}}\equiv b^{p^{n-k-s-1}} \pmod
{N_{k}}$. Let also $N_{n-s}=G$. For any $i=0,1,\dots ,n-s$ we set
$M_{i}=\bigcap_{k=i}^{n-s} N_{k}$ and claim that subgroup $M=N\cap M_{0}$ is
required. Indeed, since all members of the increasing sequence of normal
subgroups $M$, $M_{0}$, $M_{1}$, \dots, $M_{n-s-1}$, $M_{n-s}=G$ have finite
$p$-index in group $G$ we can supplement it to such sequence of normal subgroups
of $G$ all quotients of which are of order $p$. The immediately verification
shows that the resulting sequence satisfies all requirements in definition of
$(A_{1},B,\varphi ,p)$@-compatible subgroup.
\bigskip

\centerline{{\bf References}}
\medskip

\noindent\item{1.} {\it Baumslag B., Tretkoff M.} Residually  finite  $HNN$
extensions. Communic. in Algebra. 1978. Vol. 6. P.~179--194.

\noindent\item{2.} {\it Baumslag G.} On the  residual  finiteness  of
generalized  free products of nilpotent groups.  Trans. Amer. Math. Soc. 1963.
Vol. 106. P.~497--498.

 \noindent\item{3.} {\it Baumslag G., Solitar D.} Some two-generator one-relator
non-Hopfian groups. Bull. Amer. Math. Soc. 1962. Vol. 68. P.~199--201.

\noindent\item{4.} {\it Brunner A. M.} On a class of one-relator groups. Can. J.
Math. 1980. Vol.~50. P.~414--420.

\noindent\item{5.} {\it Cohen D.} Residual finiteness and  Britton's  lemma. J.
London Math. Soc.(2). 1977. Vol. 16. P.~232--234.

\noindent\item{6.}  {\it  Higman G.} Amalgams of $p$@-groups. J. Algebra. 1964.
Vol. 1. P.~301--305.

\noindent\item{7.} {\it Karrass A., Solitar D.} Subgroups of HNN groups and
groups with one defining relation. Can. J. Math. 1971. Vol. 28. P.~627--643.

\noindent\item{8.} {\it Kavutskii M.~A., Moldavanskii D.~I.} On certain class of
one-relator groups. Algebraic and discrete systems. Ivanovo. 1988. P.~35--48.
(Russian)

\noindent\item{9.} {\it Loginova E.~D.} Residual finiteness of free product of
two groups with commuting subgroups.  Sibirian Math. J. 1999. V.~40. No~2.
P.~395--407. (Russian)

\noindent\item{10.} {\it Meskin S.} Nonresidually finite one-relator groups.
Trans. Amer. Math. Soc. 1972. Vol. 164. P.~105--114.

\noindent\item{11.} {\it Moldavanskii D.~I.} Residuality a finite $p$-group of
$HNN$-extensions of finite $p$-groups. Proc. 3-rd Internat. Conf. on Algebra.
Krasnoyarsk. 1993. (Russian)

\noindent\item{12.} {\it Raptis E., Varsos D.} The residual nilpotence of
HNN-extensions with base group a finite or a f. g. abelian  group. J. of Pure
Appl. Algebra 1991. Vol.~76. P.~167--178.

\end